\magnification1200
\hfuzz=5pt
\baselineskip15pt

\def\cAq{{\overline\cA}}
\def\cDq{{\overline\cD}}

\def\[[#1]]{[\![#1]\!]}
\def\((#1)){(\!(#1)\!)}

\newread\AUX\immediate\openin\AUX=\jobname.aux
\newcount\relFnno
\def\ref#1{\expandafter\edef\csname#1\endcsname}
\ifeof\AUX\immediate\write16{\jobname.aux gibt es nicht!}\else
\input \jobname.aux
\fi\immediate\closein\AUX



\def\ignore{\bgroup
\catcode`\;=0\catcode`\^^I=14\catcode`\^^J=14\catcode`\^^M=14
\catcode`\ =14\catcode`\!=14\catcode`\"=14\catcode`\#=14\catcode`\$=14
\catcode`\&=14\catcode`\'=14\catcode`\(=14\catcode`\)=14\catcode`\*=14
\catcode`+=14\catcode`\,=14\catcode`\-=14\catcode`\.=14\catcode`\/=14
\catcode`\0=14\catcode`\1=14\catcode`\2=14\catcode`\3=14\catcode`\4=14
\catcode`\5=14\catcode`\6=14\catcode`\7=14\catcode`\8=14\catcode`\9=14
\catcode`\:=14\catcode`\<=14\catcode`\==14\catcode`\>=14\catcode`\?=14
\catcode`\@=14\catcode`\A=14\catcode`\B=14\catcode`\C=14\catcode`\D=14
\catcode`\E=14\catcode`\F=14\catcode`\G=14\catcode`\H=14\catcode`\I=14
\catcode`\J=14\catcode`\K=14\catcode`\L=14\catcode`\M=14\catcode`\N=14
\catcode`\O=14\catcode`\P=14\catcode`\Q=14\catcode`\R=14\catcode`\S=14
\catcode`\T=14\catcode`\U=14\catcode`\V=14\catcode`\W=14\catcode`\X=14
\catcode`\Y=14\catcode`\Z=14\catcode`\[=14\catcode`\\=14\catcode`\]=14
\catcode`\^=14\catcode`\_=14\catcode`\`=14\catcode`\a=14\catcode`\b=14
\catcode`\c=14\catcode`\d=14\catcode`\e=14\catcode`\f=14\catcode`\g=14
\catcode`\h=14\catcode`\i=14\catcode`\j=14\catcode`\k=14\catcode`\l=14
\catcode`\m=14\catcode`\n=14\catcode`\o=14\catcode`\p=14\catcode`\q=14
\catcode`\r=14\catcode`\s=14\catcode`\t=14\catcode`\u=14\catcode`\v=14
\catcode`\w=14\catcode`\x=14\catcode`\y=14\catcode`\z=14\catcode`\{=14
\catcode`\|=14\catcode`\}=14\catcode`\~=14\catcode`\^^?=14
\Ignoriere}
\def\Ignoriere#1\;{\egroup}

\newcount\itemcount
\def\resetitem{\global\itemcount0}\resetitem
\newcount\Itemcount
\Itemcount0
\newcount\maxItemcount
\maxItemcount=0

\def\FILTER\fam\itfam\tenit#1){#1}

\def\Item#1{\global\advance\itemcount1
\edef\TEXT{{\it\romannumeral\itemcount)}}%
\ifx?#1?\relax\else
\ifnum#1>\maxItemcount\global\maxItemcount=#1\fi
\expandafter\ifx\csname I#1\endcsname\relax\else
\edef\testA{\csname I#1\endcsname}
\expandafter\expandafter\def\expandafter\testB\testA
\edef\testC{\expandafter\FILTER\testB}
\edef\testD{\csname0\testC0\endcsname}\fi
\edef\testE{\csname0\romannumeral\itemcount0\endcsname}
\ifx\testD\testE\relax\else
\immediate\write16{I#1 hat sich geaendert!}\fi
\expandwrite\AUX{\neverexpand\ref{I#1}{\TEXT}}\fi
\item{\TEXT}}

\def\today{\number\day.~\ifcase\month\or
  Januar\or Februar\or M{\"a}rz\or April\or Mai\or Juni\or
  Juli\or August\or September\or Oktober\or November\or Dezember\fi
  \space\number\year}
\font\sevenex=cmex7
\scriptfont3=\sevenex
\font\fiveex=cmex10 scaled 500
\scriptscriptfont3=\fiveex
\def\A{{\bf A}}
\def\G{{\bf G}}
\def\P{{\bf P}}

\def\phi{\varphi}
\def\epsilon{\varepsilon}

\def\theta{\vartheta}
\font\lams=lams3
\def\uinto{\lower1.7pt\hbox{%
\vbox{\offinterlineskip
\hbox{\lams\char"7A}%
\hbox{\vbox to 7.5pt{\leaders\vrule width0.2pt\vfill}%
\kern-.3pt\hbox{\lams\char"76}}}}}
\def\uauf{\lower1.7pt\hbox to 3pt{%
\vbox{\offinterlineskip
\hbox{\vbox to 8.5pt{\leaders\vrule width0.2pt\vfill}%
\kern-.3pt\hbox{\lams\char"76}\kern-0.3pt%
$\raise1pt\hbox{\lams\char"76}$}}\hfil}}

\def\title#1{\par
{\baselineskip1.5\baselineskip\rightskip0pt plus 5truecm
\leavevmode\vskip0truecm\noindent\font\BF=cmbx10 scaled \magstep2\BF #1\par}
\vskip1truecm
\leftline{\font\CSC=cmcsc10
{\CSC Friedrich Knop}}
\leftline{Department of Mathematics, Rutgers University, New Brunswick NJ
08903, USA}
\leftline{knop@math.rutgers.edu}
\vskip1truecm
\par}

\def\cite#1{\expandafter\ifx\csname#1\endcsname\relax
{\bf?}\immediate\write16{#1 ist nicht definiert!}\else\csname#1\endcsname\fi}
\def\expandwrite#1#2{\edef\next{\write#1{#2}}\next}
\def\neverexpand{\noexpand\noexpand\noexpand}
\def\strip#1\ {}
\def\ncite#1{\expandafter\ifx\csname#1\endcsname\relax
{\bf?}\immediate\write16{#1 ist nicht definiert!}\else
\expandafter\expandafter\expandafter\strip\csname#1\endcsname\fi}
\newwrite\AUX
\immediate\openout\AUX=\jobname.aux
\font\eightrm=cmr8\font\sixrm=cmr6
\font\eighti=cmmi8
\font\eightit=cmti8
\font\eightbf=cmbx8
\font\eightcsc=cmcsc10 scaled 833
\def\eightpoint{%
\textfont0=\eightrm\scriptfont0=\sixrm\def\rm{\fam0\eightrm}%
\textfont1=\eighti
\textfont\bffam=\eightbf\def\bf{\fam\bffam\eightbf}%
\textfont\itfam=\eightit\def\it{\fam\itfam\eightit}%
\def\csc{\eightcsc}%
\setbox\strutbox=\hbox{\vrule height7pt depth2pt width0pt}%
\normalbaselineskip=0,8\normalbaselineskip\normalbaselines\rm}
\newcount\absFnno\absFnno1
\write\AUX{\relFnno1}
\newif\ifMARKE\MARKEtrue
{\catcode`\@=11
\gdef\footnote{\ifMARKE\edef\@sf{\spacefactor\the\spacefactor}\/%
$^{\cite{Fn\the\absFnno}}$\@sf\fi
\MARKEtrue
\insert\footins\bgroup\eightpoint
\interlinepenalty100\let\par=\endgraf
\leftskip=0pt\rightskip=0pt
\splittopskip=10pt plus 1pt minus 1pt \floatingpenalty=20000\smallskip
\item{$^{\cite{Fn\the\absFnno}}$}%
\expandwrite\AUX{\neverexpand\ref{Fn\the\absFnno}{\neverexpand\the\relFnno}}%
\global\advance\absFnno1\write\AUX{\advance\relFnno1}%
\bgroup\strut\aftergroup\@foot\let\next}}
\skip\footins=12pt plus 2pt minus 4pt
\dimen\footins=30pc
\output={\plainoutput\immediate\write\AUX{\relFnno1}}
\newcount\Abschnitt\Abschnitt0
\def\beginsection#1. #2 \par{\advance\Abschnitt1%
\vskip0pt plus.10\vsize\penalty-250
\vskip0pt plus-.10\vsize\bigskip\vskip\parskip
\edef\TEST{\number\Abschnitt}
\expandafter\ifx\csname#1\endcsname\TEST\relax\else
\immediate\write16{#1 hat sich geaendert!}\fi
\expandwrite\AUX{\neverexpand\ref{#1}{\TEST}}
\leftline{\marginnote{#1}\bf\number\Abschnitt. \ignorespaces#2}%
\nobreak\smallskip\noindent\SATZ1\GNo0}
\def\Proof:{\par\noindent{\it Proof:}}
\def\Remark:{\ifdim\lastskip<\medskipamount\removelastskip\medskip\fi
\noindent{\bf Remark:}}
\def\Remarks:{\ifdim\lastskip<\medskipamount\removelastskip\medskip\fi
\noindent{\bf Remarks:}}
\def\Definition:{\ifdim\lastskip<\medskipamount\removelastskip\medskip\fi
\noindent{\bf Definition:}}
\def\Example:{\ifdim\lastskip<\medskipamount\removelastskip\medskip\fi
\noindent{\bf Example:}}
\def\Examples:{\ifdim\lastskip<\medskipamount\removelastskip\medskip\fi
\noindent{\bf Examples:}}
\newif\ifmarginalnotes\marginalnotesfalse
\newif\ifmarginalwarnings\marginalwarningstrue

\def\marginnote#1{\ifmarginalnotes\hbox to 0pt{\eightpoint\hss #1\ }\fi}

\def\strutdepth{\dp\strutbox}
\def\Randbem#1#2{\ifmarginalwarnings
{#1}\strut
\setbox0=\vtop{\eightpoint
\rightskip=0pt plus 6mm\hfuzz=3pt\hsize=16mm\noindent\leavevmode#2}%
\vadjust{\kern-\strutdepth
\vtop to \strutdepth{\kern-\ht0
\hbox to \hsize{\kern-16mm\kern-6pt\box0\kern6pt\hfill}\vss}}\fi}

\def\Zitat!{\Randbem{\bf?}{\bf Zitat}}

\newcount\SATZ\SATZ1
\def\proclaim #1. #2\par{\ifdim\lastskip<\medskipamount\removelastskip
\medskip\fi
\noindent{\bf#1.\ }{\it#2}\Par
\ifdim\lastskip<\medskipamount\removelastskip\goodbreak\medskip\fi}
\def\Aussage#1{\expandafter\def\csname#1\endcsname##1.{\resetitem
\ifx?##1?\relax\else
\edef\TEST{#1\penalty10000\ \number\Abschnitt.\number\SATZ}
\expandafter\ifx\csname##1\endcsname\TEST\relax\else
\immediate\write16{##1 hat sich geaendert!}\fi
\expandwrite\AUX{\neverexpand\ref{##1}{\TEST}}\fi
\proclaim {\marginnote{##1}\number\Abschnitt.\number\SATZ. #1\global\advance\SATZ1}.}}
\Aussage{Theorem}
\Aussage{Proposition}
\Aussage{Corollary}
\Aussage{Lemma}
\font\la=lasy10
\def\strich{\hbox{$\vcenter{\hbox
to 1pt{\leaders\hrule height -0,2pt depth 0,6pt\hfil}}$}}
\def\dashedrightarrow{\hbox{%
\hbox to 0,5cm{\leaders\hbox to 2pt{\hfil\strich\hfil}\hfil}%
\kern-2pt\hbox{\la\char\string"29}}}

\def\Bindestrich{\penalty10000-\hskip0pt}
\let\_=\Bindestrich
\def\.{{\sfcode`.=1000.}}
\def\Links#1{\llap{$\scriptstyle#1$}}

\def\Par{\par}
\def\:={\mathrel{\raise0,9pt\hbox{.}\kern-2,77779pt
\raise3pt\hbox{.}\kern-2,5pt=}}
\def\=:{\mathrel{=\kern-2,5pt\raise0,9pt\hbox{.}\kern-2,77779pt
\raise3pt\hbox{.}}} \def\mod{/\mskip-5mu/}
\def\into{\hookrightarrow}
\def\pfeil{\rightarrow}

\def\Pf#1{\buildrel#1\over\longrightarrow}

\def\Ugleich{\hbox{$\cup$\kern.5pt\vrule depth -0.5pt}}
\def\|#1|{\mathop{\rm#1}\nolimits}
\def\<{\langle}
\def\>{\rangle}
\let\Times=\times
\def\times{\mathop{\Times}}
\let\Otimes=\otimes
\def\otimes{\mathop{\Otimes}}
\catcode`\@=11
\def\hex#1{\ifcase#1 0\or1\or2\or3\or4\or5\or6\or7\or8\or9\or A\or B\or
C\or D\or E\or F\else\message{Warnung: Setze hex#1=0}0\fi}
\def\fontdef#1:#2,#3,#4.{%
\alloc@8\fam\chardef\sixt@@n\FAM
\ifx!#2!\else\expandafter\font\csname text#1\endcsname=#2
\textfont\the\FAM=\csname text#1\endcsname\fi
\ifx!#3!\else\expandafter\font\csname script#1\endcsname=#3
\scriptfont\the\FAM=\csname script#1\endcsname\fi
\ifx!#4!\else\expandafter\font\csname scriptscript#1\endcsname=#4
\scriptscriptfont\the\FAM=\csname scriptscript#1\endcsname\fi
\expandafter\edef\csname #1\endcsname{\fam\the\FAM\csname text#1\endcsname}
\expandafter\edef\csname hex#1fam\endcsname{\hex\FAM}}
\catcode`\@=12 

\fontdef Ss:cmss10,,.
\fontdef Fr:eufm10,eufm7,eufm5.


\def\fm{{\Fr m}}

\fontdef bbb:msbm10,msbm7,msbm5.
\fontdef mbf:cmmib10,cmmib7,.
\fontdef msa:msam10,msam7,msam5.
\def\CC{{\bbb C}}

\def\NN{{\bbb N}}

\def\ZZ{{\bbb Z}}
\def\cA{{\cal A}}\def\cD{{\cal D}}

\def\cO{{\cal O}}\def\cP{{\cal P}}
\def\cT{{\cal T}}

\mathchardef\leer=\string"0\hexbbbfam3F
\mathchardef\subsetneq=\string"3\hexbbbfam24
\mathchardef\semidir=\string"2\hexbbbfam6E
\mathchardef\dirsemi=\string"2\hexbbbfam6F
\mathchardef\haken=\string"2\hexmsafam78
\mathchardef\auf=\string"3\hexmsafam10
\let\OL=\overline
\def\overline#1{{\hskip1pt\OL{\hskip-1pt#1\hskip-.3pt}\hskip.3pt}}


\def\Gq{{\overline{G}}}

%
\newdimen\Parindent
\Parindent=\parindent


\abovedisplayskip 9.0pt plus 3.0pt minus 3.0pt
\belowdisplayskip 9.0pt plus 3.0pt minus 3.0pt
\newdimen\Grenze\Grenze2\Parindent\advance\Grenze1em
\newdimen\Breite
\newbox\DpBox
\def\NewDisplay#1
#2$${\Breite\hsize\advance\Breite-\hangindent
\setbox\DpBox=\hbox{\hskip2\Parindent$\displaystyle{%
\ifx0#1\relax\else\eqno{#1}\fi#2}$}%
\ifnum\predisplaysize<\Grenze\abovedisplayskip\abovedisplayshortskip
\belowdisplayskip\belowdisplayshortskip\fi
\global\futurelet\nexttok\WEITER}
\def\WEITER{\ifx\nexttok\qed\expandafter\leftQEDdisplay
\else\leftdisplay\fi}
\def\leftdisplay{\hskip-\hangindent\leftline{\box\DpBox}$$}
\def\leftQEDdisplay{\hskip-\hangindent
\line{\copy\DpBox\hfill\lower\dp\DpBox\copy\QEDbox}%
\belowdisplayskip0pt$$\bigskip\let\nexttok=}
\everydisplay{\NewDisplay}
\newcount\GNo\GNo=0
\newcount\maxEqNo\maxEqNo=0
\def\eqno#1{%
\global\advance\GNo1
\edef\FTEST{$(\number\Abschnitt.\number\GNo)$}
\ifx?#1?\relax\else
\ifnum#1>\maxEqNo\global\maxEqNo=#1\fi%
\expandafter\ifx\csname E#1\endcsname\FTEST\relax\else
\immediate\write16{E#1 hat sich geaendert!}\fi
\expandwrite\AUX{\neverexpand\ref{E#1}{\FTEST}}\fi
\llap{\hbox to 40pt{\marginnote{#1}\FTEST\hfill}}}

\catcode`@=11
\def\eqalignno#1{\null\!\!\vcenter{\openup\jot\m@th\ialign{\eqno{##}\hfil
&\strut\hfil$\displaystyle{##}$&$\displaystyle{{}##}$\hfil\crcr#1\crcr}}\,}
\catcode`@=12

\newbox\QEDbox
\newbox\nichts\setbox\nichts=\vbox{}\wd\nichts=2mm\ht\nichts=2mm
\setbox\QEDbox=\hbox{\vrule\vbox{\hrule\copy\nichts\hrule}\vrule}
\def\qed{\leavevmode\unskip\hfil\null\nobreak\hfill\copy\QEDbox\medbreak}
\newdimen\HIindent
\newbox\HIbox
\def\setHI#1{\setbox\HIbox=\hbox{#1}\HIindent=\wd\HIbox}
\def\HI#1{\par\hangindent\HIindent\hangafter=0\noindent\leavevmode
\llap{\hbox to\HIindent{#1\hfil}}\ignorespaces}

\newdimen\maxSpalbr
\newdimen\altSpalbr
\newcount\Zaehler


\newif\ifxxx

{\catcode`/=\active

\gdef\beginrefs{%
\xxxfalse
\catcode`/=\active
\def/{\string/\ifxxx\hskip0pt\fi}
\def\TText##1{{\xxxtrue\tt##1}}
\expandafter\ifx\csname Spaltenbreite\endcsname\relax
\def\Spaltenbreite{1cm}\immediate\write16{Spaltenbreite undefiniert!}\fi
\expandafter\altSpalbr\Spaltenbreite
\maxSpalbr0pt
\gdef\alt{}
\def\\##1\relax{%
\gdef\neu{##1}\ifx\alt\neu\global\advance\Zaehler1\else
\xdef\alt{\neu}\global\Zaehler=1\fi\xdef\SigText{##1\the\Zaehler}}
\def\L|Abk:##1|Sig:##2|Au:##3|Tit:##4|Zs:##5|Bd:##6|S:##7|J:##8|xxx:##9||{%
\def\SigText{##2}\global\setbox0=\hbox{##2\relax}
\edef\TEST{[\SigText]}
\expandafter\ifx\csname##1\endcsname\TEST\relax\else
\immediate\write16{##1 hat sich geaendert!}\fi
\expandwrite\AUX{\neverexpand\ref{##1}{\TEST}}
\setHI{[\SigText]\ }
\ifnum\HIindent>\maxSpalbr\maxSpalbr\HIindent\fi
\ifnum\HIindent<\altSpalbr\HIindent\altSpalbr\fi
\HI{\marginnote{##1}[\SigText]}
\ifx-##3\relax\else{##3}: \fi
\ifx-##4\relax\else{##4}{\sfcode`.=3000.} \fi
\ifx-##5\relax\else{\it ##5\/} \fi
\ifx-##6\relax\else{\bf ##6} \fi
\ifx-##8\relax\else({##8})\fi
\ifx-##7\relax\else, {##7}\fi
\ifx-##9\relax\else, \TText{##9}\fi\Par}
\def\B|Abk:##1|Sig:##2|Au:##3|Tit:##4|Reihe:##5|Verlag:##6|Ort:##7|J:##8|xxx:##9||{%
\def\SigText{##2}\global\setbox0=\hbox{##2\relax}
\edef\TEST{[\SigText]}
\expandafter\ifx\csname##1\endcsname\TEST\relax\else
\immediate\write16{##1 hat sich geaendert!}\fi
\expandwrite\AUX{\neverexpand\ref{##1}{\TEST}}
\setHI{[\SigText]\ }
\ifnum\HIindent>\maxSpalbr\maxSpalbr\HIindent\fi
\ifnum\HIindent<\altSpalbr\HIindent\altSpalbr\fi
\HI{\marginnote{##1}[\SigText]}
\ifx-##3\relax\else{##3}: \fi
\ifx-##4\relax\else{##4}{\sfcode`.=3000.} \fi
\ifx-##5\relax\else{(##5)} \fi
\ifx-##7\relax\else{##7:} \fi
\ifx-##6\relax\else{##6}\fi
\ifx-##8\relax\else{ ##8}\fi
\ifx-##9\relax\else, \TText{##9}\fi\Par}
\def\Pr|Abk:##1|Sig:##2|Au:##3|Artikel:##4|Titel:##5|Hgr:##6|Reihe:{%
\def\SigText{##2}\global\setbox0=\hbox{##2\relax}
\edef\TEST{[\SigText]}
\expandafter\ifx\csname##1\endcsname\TEST\relax\else
\immediate\write16{##1 hat sich geaendert!}\fi
\expandwrite\AUX{\neverexpand\ref{##1}{\TEST}}
\setHI{[\SigText]\ }
\ifnum\HIindent>\maxSpalbr\maxSpalbr\HIindent\fi
\ifnum\HIindent<\altSpalbr\HIindent\altSpalbr\fi
\HI{\marginnote{##1}[\SigText]}
\ifx-##3\relax\else{##3}: \fi
\ifx-##4\relax\else{##4}{\sfcode`.=3000.} \fi
\ifx-##5\relax\else{In: \it ##5}. \fi
\ifx-##6\relax\else{(##6)} \fi\PrII}
\def\PrII##1|Bd:##2|Verlag:##3|Ort:##4|S:##5|J:##6|xxx:##7||{%
\ifx-##1\relax\else{##1} \fi
\ifx-##2\relax\else{\bf ##2}, \fi
\ifx-##4\relax\else{##4:} \fi
\ifx-##3\relax\else{##3} \fi
\ifx-##6\relax\else{##6}\fi
\ifx-##5\relax\else{, ##5}\fi
\ifx-##7\relax\else, \TText{##7}\fi\Par}
\bgroup
\baselineskip12pt
\parskip2.5pt plus 1pt
\hyphenation{Hei-del-berg Sprin-ger}
\sfcode`.=1000
\beginsection References. References

}}

\def\endrefs{%
\expandwrite\AUX{\neverexpand\ref{Spaltenbreite}{\the\maxSpalbr}}
\ifnum\maxSpalbr=\altSpalbr\relax\else
\immediate\write16{Spaltenbreite hat sich geaendert!}\fi
\egroup\write16{Letzte Gleichung: E\the\maxEqNo}
\write16{Letzte Aufzaehlung: I\the\maxItemcount}}



\title{Graded cofinite rings of differential operators}

{\narrower\noindent We classify subalgebras of a ring of differential
operators which are big in the following sense: the extension of
associated graded rings is finite. We show that these subalgebras
correspond, up to automorphism, to uniformly ramified finite
morphisms. This generalizes a theorem of Levasseur-Stafford on the
generators of the invariants of a Weyl algebra under a finite group.

}

\beginsection Intro. Introduction

In this paper we study subalgebras $\cA$ of the algebra $\cD(X)$ of
differential operators on a smooth variety $X$ which are big in the
following sense: using the order of a differential operator, the ring
$\cD(X)$ is equipped with a filtration. Its associated graded algebra
$\cDq(X)$ is commutative and can be regarded as the set of regular
functions on the cotangent bundle of $X$. The subalgebra $\cA$
inherits a filtration from $\cD(X)$ and its associated graded algebra
$\cAq$ is a subalgebra of $\cDq(X)$. We call $\cA$ {\it graded
cofinite} in $\cD(X)$ if $\cDq(X)$ is a finitely generated
$\cAq$\_module.

Our guiding example of a graded cofinite subalgebra is the
algebra of invariants $\cD(X)^W$ where $W$ is a finite group acting on
$X$.

Other examples can be constructed as follows. Let $\phi:X\pfeil Y$ be
a finite dominant morphism onto a normal variety $Y$. Then we put
$$
\cD(X,Y)=\{D\in\cD(X)\mid D(\cO(Y))\subseteq\cO(Y)\}.
$$
We show (\cite{gradedcofinite}) that this subalgebra is graded
cofinite if and only if the ramification of $\phi$ is uniform, i.e.,
the ramification degree of $\phi$ along a divisor $D\subset X$ depends
only on the image $\phi(D)$.

It should be noted that these two constructions are in fact more or
less equivalent. In \cite{finite} we show that
$\cD(X)^W=\cD(X,X/W)$. Conversely, we show in \cite{reduction} that
$\cD(X,Y)=\cD(\tilde X)^W$ where $\tilde X\pfeil X$ is a suitable
finite cover of $X$ and $W$ is a finite group acting on $\tilde X$.

Our main result is that up to automorphisms every graded cofinite
subalgebra is of form above:

\Theorem. Let $X$ be a smooth variety and $\cA$ a graded cofinite
subalgebra of $\cD(X)$. Then there is an automorphism $\Phi$ of
$\cD(X)$, inducing the identity on $\cDq(X)$, such that
$\cA=\Phi\cD(X,Y)$ for some uniformly ramified morphism $\phi:X\pfeil
Y$.

The main motivation for this notion came from the following result of
Levasseur and Stafford: let $W$ be a finite group acting linearly on a
vector space $V$. Then $\cD(V)^W$ is generated by the $W$\_invariant
functions $\cO(V)^W$ and the $W$\_invariant constant coefficient
differential operators $S^*(V)^W$. For general varieties $X$, there is
no notion of constant coefficient differential operators. Since the
algebra generated by $\cO(V)^W$ and $S^*(V)^W$ is clearly graded
cofinite our main theorem can be seen as a non\_linear generalization
of the theorem of Levasseur\_Stafford.

Our main theorem has several application concerning generating
elements of rings of $W$\_invariant differential operators which go
beyond the theorem of Levasseur\_Stafford. For example, we prove that
$\cD(X)^W$ can be generated by at most $2n+1$ elements when $V$ is an
$n$\_dimensional representation of $W$. Moreover, we establish a kind of
Galois correspondence for graded cofinite subalgebras. Finally, we
determine all graded cofinite subalgebras of $\cD(\A^1)$, the Weyl
algebra in two generators.

The proof consists essentially of five steps: 1. We show the
aforementioned claim that $\cD(X,Y)$ is graded cofinite if and only if
$\phi$ is uniformly ramified. 2. Then we show that under these
conditions $\cD(X,Y)$ is a simple ring. Here we follow an argument of
Wallach~\cite{Wa}. 3. We show that the theorem holds over the generic
point of $X$. 4. Then we construct the automorphism
$\Phi$. This is the most tedious part of the paper and rests on
explicit computations in codimension one. 5. Finally, we paste all
this information together by showing that two graded cofinite
subalgebras $\cA\subseteq\cA'$ which coincide generically and for
which $\cA'$ is a simple ring are actually equal. Here we follow the
argument in \cite{LS}.

Finally, it should be mentioned that the actual main \cite{main} is more
general in that it allows for certain singularities of $X$.

\medskip \noindent{\bf Acknowledgment:} This work started while the
author was guest of the CRM, Montr\'eal, in Summer 1997 and continued
during a stay at the University of Freiburg in 2004. The author thanks
both institutions for their hospitality. Last not least, the author
would like to thank the referee for an excellent job. In particular,
the shorter proof of \cite{finite} was pointed out by him/her.

\beginsection basechange. Graded cofinite subalgebras: definition and
base change

All varieties and algebras will be defined over $\CC$. Moreover, varieties are
irreducible by definition.

Recall that a $\CC$\_linear endomorphism $D$ of a commutative algebra $B$ is a
{\it differential operator of order $\le d$} if
$$
[b_0,[b_1,\ldots[b_d,D]\ldots]]=0\quad\hbox{for all}\quad
 b_0,b_1,\ldots,b_d\in B.
$$
Let $\cD(B)_{\le d}$ be the set of differential operators of order
$\le d$ and and $\cD(B)=\bigcup_d\cD(B)_{\le d}$. Then $\cD(B)$ is a
filtered algebra, i.e., $\cD(B)_{\le d}\cD(B)_{\le
e}\subseteq\cD(B)_{\le d+e}$ for all integers $d$ and $e$. Let
$\cDq(B)$ be its {\it associated graded algebra}, i.e.,
$\cDq(X):=\oplus_d\cDq(X)_d$ with $\cDq(X)_d=\cD(X)_{\le
d}/\cD(X)_{\le d-1}$. This is a graded {\it commutative} algebra.
If $X$ is a variety with ring of functions $\cO(X)$ then we define
$\cD(X)=\cD(\cO(X))$.

Every subalgebra $\cA\subseteq\cD(X)$ inherits the filtration by
$\cA_{\le d}=\cA\cap\cD(X)_{\le d}$. This way, the associated graded
algebra $\cAq$ is a subalgebra of $\cDq(X)$ and we define:

\Definition: A subalgebra $\cA$ of $\cD(X)$ is called {\it graded
cofinite} if $\cDq(X)$ is a finitely generated $\cAq$\_module.

\Example: Let $W$ be a finite group acting on $X$ and assume $\cDq(X)$
to be finitely generated (e.g. $X$ smooth). We claim that
$\cA=\cD(X)^W$ is graded cofinite in $\cD(X)$. In fact, since $W$ is
linearly reductive, we have $\cAq=\cDq(X)^W$ which is well known to be
cofinite in $\cDq(X)$.

\medskip
The ring $A:=\cA_{\le0}=\cA\cap\cO(X)$ is called the {\it base} of
$\cA$.

\Proposition finitebase. Let $\cA\subseteq\cD(X)$ be graded
cofinite. Then the base $A$ of $\cA$ is a finitely generated algebra
which is cofinite in $\cO(X)$. In other words, if $Y=\|Spec|A$ then
$X\pfeil Y$ is a finite surjective morphism of affine varieties.

\Proof: Since $\cAq$ is cofinite in $\cDq(X)$, its 0-component $A$ is
cofinite in the 0-component $\cO(X)$ of $\cDq(X)$. Now the assertion
follows from the following lemma.\qed

\Lemma finitecomm. Let $A\subseteq B$ be an integral extension of
commutative $\CC$\_algebras. Assume $B$ is a finitely generated
algebra. Then $A$ is finitely generated as well and $A$ is cofinite in
$B$.

\Proof: This is the Artin-Tate lemma. For a proof see \cite{Eisen}
p.~143.\qed

In the sequel we need some auxiliary results concerning the behavior of
$\cA$ with respect to extension of scalars. Let $X$ be an affine
variety, $B:=\cO(X)$, and $J\subseteq B$ an ideal. Let $\hat B$ be the
$J$\_adic completion of $B$ and $\hat X:=\|Spec|\hat B$. Let
$\cD_c(\hat X)\subseteq\|End|_\CC^{\rm cont}(\hat B)$ be the algebra
of {\it continuous} differential operators on $\hat X$. We show that
this is also the algebra of differential operators on $X$ with
coefficients in $\hat B$. More precisely:

\Lemma completion0. Fix $d\ge0$. Then the left $J$\_adic topology and the right
$J$\_adic topology of $\cD(X)_{\le d}$ coincide. Its completion with
respect to this topology equals $\cD_c(\hat X)_{\le d}$. In particular,
the two natural maps
$$
\hat B\otimes_B\cD(X)\rightarrow\cD_c(\hat X)\quad{\rm and}\quad
\cD(X)\otimes_B\hat B\pfeil\cD_c(\hat X)
$$
are isomorphisms of filtered vector spaces.

\Proof: We recall Grothendieck's description of $\cD(X)$: let $\delta$ be
the kernel of the multiplication map $B\otimes_\CC B\pfeil B$. It is
the ideal of $C:=B\otimes B$ generated by all elements of the form
$b\otimes1-1\otimes b$, $b\in B$. Let $\cP^d:=C/\delta^{d+1}$. This is
a $C$\_module, i.e., carries a left and a right $B$\_module
structure. Moreover, it is a finitely generated module with respect to
both structures. Now we have $\cD(X)_{\le d}=\|Hom|_B(\cP^d,B)$ where
we use the left $B$\_module structure of $\cP^d$.

Now consider the completed ring $\hat B$. Then $\|End|_\CC(\hat
B)=\|Hom|_{\hat B}(\hat B\otimes_\CC\hat B,\hat B)$. It is easy to see
that the continuous endomorphisms correspond exactly to those
homomorphisms $\hat B\otimes_\CC\hat B\pfeil\hat B$ which extend to
the completed tensor product $\hat C:=\hat B\hat\Otimes_\CC\hat
B$. Thus, $\|End|^{\rm cont}_\CC(\hat B)=\|Hom|_{\hat B}(\hat C,\hat
B)$. Let $\hat\delta$ be the kernel of $\hat C=\hat
B\hat\Otimes_\CC\hat B\pfeil\hat B$ and $\hat\cP^d=\hat
C/\hat\delta^{d+1}$. Then $\cD_c(\hat X)_{\le d}=\|Hom|_{\hat
B}(\hat\cP^d,\hat B)$.

Let $K:=J\otimes B+B\otimes J\subseteq C$. Then $\hat C$ is the
$K$\_adic completion of $C$. Moreover, $\hat\delta$ is the $K$\_adic
completion of $\delta$. Thus, everything boils down to the
following statement: the left $J$\_adic, the right $J$\_adic, and the
$K$\_adic topologies of $\cP^d$ all coincide.

For $b\in B$ we have $1\otimes b=b\otimes 1+c$ with
$c=b\otimes1-1\otimes b\in\delta$. Thus $B\otimes J\subseteq
J\otimes B+\delta$ and, for any $n\ge d$,
$$
J^n\otimes B\subseteq K^n\subseteq(J\otimes B+\delta)^n\subseteq
J^{n-d}\otimes B+\delta^{d+1}.
$$
This shows that the left $J$\_adic and the $K$\_adic topologies of
$\cP^d$ coincide. The argument for the right $J$\_adic topology is the
same.\qed

Now let $\cA\subseteq\cD(X)$ be graded cofinite with base $A$. Let
$I\subseteq A$ be an ideal and let $\hat A$ be the $I$\_adic
completion of $A$. Set $J:=IB\subseteq B$. Since $J^n=I^nB$, the
$I$\_adic completion $\hat B$ of $B$ is the same as its $J$\_adic
completion.

\Corollary completion. Let $\hat\cA\subseteq\cD_c(\hat X)$ be the
subalgebra generated by $\cA$ and $\hat A$. Then $\hat\cA$ is a graded
cofinite subalgebra of $\cD_c(\hat X)$ with base $\hat A$. Moreover,
the maps
$$5
\hat A\otimes_A\cA\pfeil\hat\cA\quad{\rm and}\quad
\cA\otimes_A\hat A\pfeil\hat\cA
$$
are isomorphisms of filtered vector spaces.

\Proof: Redefine $\hat\cA$ to be the closure of $\cA$ in $\cD_c(\hat
X)$ with respect to either left or right $J$\_adic topology. Then
\cite{completion0} implies that the maps \cite{E5} are isomorphisms. In
particular, $\hat\cA$ is an algebra and therefore the algebra
generated by $\hat A$ and $\cA$.\qed

\noindent
Now we deduce the same thing for \'etale morphisms. Again,
let $\cA\subseteq\cD(X)$ be graded cofinite with base $A$. Let $\tilde
Y\pfeil Y:=\|Spec|A$ be an \'etale morphism where $\tilde Y$ is
another affine variety. Then also $\tilde X:=\tilde Y\times_YX\pfeil
X$ is \'etale. Now put $B:=\cO(X)$, $\tilde A:=\cO(\tilde Y)$, and
$\tilde B:=\cO(\tilde X)=\tilde A\otimes_AB$. Then
$$7
\tilde B\otimes_B\cD(X)\Pf\sim\cD(\tilde X)\quad{\rm and}\quad
\cD(X)\otimes_B\tilde B\Pf\sim\cD(\tilde X)
$$
are isomorphisms. For a proof see \cite{Masson} Thm.~2.2.10,
Prop.~2.2.12 or \cite{Sch2}~Thm.4.2. Both references state that the
first isomorphism is an isomorphism of filtered rings, i.e., that
there is an isomorphism on the associated graded level. This, in turn,
implies the second isomorphism.

\Lemma etale. Let $\cA$, $A$, $\tilde A$, and $\tilde X$ be as above.
Let $\tilde\cA\subseteq\cD(\tilde X)$ be the subalgebra
generated by $\cA$ and $\tilde A$. Then $\tilde\cA$ is a graded
cofinite subalgebra of $\cD(\tilde X)$ with base $\tilde A$. Moreover,
the maps
$$
\tilde A\otimes_A\cA\pfeil\tilde\cA\quad{\rm and}\quad
\cA\otimes_A\tilde A\pfeil\tilde\cA
$$
are isomorphisms of filtered vector spaces.

\Proof: We start with a general remark. Let $Z$ be an affine
variety. Let $\fm_z\subset\cO(Z)$ be the maximal ideal corresponding
to a point $z\in Z$. It is known that $\fm_z$\_adic completion is
exact on finitely generated $\cO(Z)$\_modules. Moreover, a finitely
generated $\cO(Z)$\_module $M$ is $0$ if and only if it is so after
$\fm_z$\_adic completion for every $z\in Z$. Now let $N\subseteq M$ be
a submodule, $N'$ another $\cO(Z)$\_module and $N'\pfeil M$ an
$\cO(Z)$\_homomorphism. Then one sees from the remarks above that
$\phi$ induces an isomorphism of $N'$ onto $N$ if and only if this is
so after $\fm_z$\_adic completion for every $z\in Z$.

We apply this to $Z=\tilde X$ and $M=\cD(\tilde X)_{\le d}$. Let
$\tilde x\in\tilde X$ with image $x\in X$. Since $\tilde X\pfeil X$ is
\'etale, the $\fm_{\tilde x}$\_adic completion of $\cD(\tilde X)_{\le d}$ is
the same as the $\fm_x$\_adic completion of $\cD(X)_{\le d}$. Thus,
\cite{completion} implies that the two homomorphisms
$$
\tilde A\otimes_A\cA_{\le d}\pfeil\cD(\tilde X)_{\le d}\quad{\rm and}\quad
\cA_{\le d}\otimes_A\tilde A\pfeil\cD(\tilde X)_{\le d}
$$
are injective with the same image after $\fm_{\tilde x}$\_adic
completion for every $\tilde x\in\tilde X$.\qed

\Corollary. Let $\cA\subseteq\cD(X)$ be graded cofinite with base
$A$ and let $S\subseteq A$ be a multiplicatively closed subset defining
localizations $A_S\subseteq B_S$ (with $B:=\cO(X)$). Let
$\cA_S\subseteq\cD(B_S)$ be the subalgebra generated by $\cA$ and
$A_S$. Then $\cA_S$ is a graded cofinite subalgebra of $\cD(B_S)$ with
base $A_S$. Moreover, the maps $A_S\otimes_A\cA\pfeil\cA_S$ and
$\cA\otimes_A A_S\pfeil\cA_S$ are filtered isomorphisms.

\Proof: If $S$ is finite then $A\pfeil A_S$ is an open embedding, in
particular \'etale. It follows from \cite{etale} that $\cA_S$ has base
$A_S$. For the general case use that $S$ is the union of its finite
subsets and that all objects behave well under inductive limits.\qed

\noindent An important consequence is that we can ``normalize'' graded
cofinite subalgebras.

\Corollary normalization. Let $X$ be normal and $\cA\subseteq\cD(X)$
be a graded cofinite subalgebra. Let $A'$ be the normalization of the
base $A$, regarded as a subalgebra of $\cO(X)$. Let
$\cA'\subseteq\cD(X)$ be the subalgebra generated by $\cA$ and
$A'$. Then $\cA'$ is a graded cofinite subalgebra of $\cD(X)$ with
base $A'$.

\Proof: Let $B:=\cO(X)$. Both algebras $A$ and $A'$ have the same
quotient field $K=A_S$ with $S=A\setminus\{0\}$.  Thus we have
$A'\subseteq\cA'\cap B\subseteq\cA_K\cap B=K\cap B=A'$.\qed

\Remark: It is possible to combine \'etale base change, localization,
and completion. More precisely, we will use this twice in the
following situation: let $\cA\subseteq\cD(X)$ be graded cofinite with
base $A$ and assume $X$ and $Y=\|Spec|A$ to be normal. Let $\tilde
Y\pfeil Y$ be \'etale and $\tilde D\subseteq\tilde Y$ a prime
divisor. Take $\hat A$ to be the completion of the local ring
$\cO_{\tilde Y,\tilde D}$. Then $\hat A\cong E\[[t]]$ is a discrete
valuation ring with $E=\CC(\tilde D)$ and $\hat B=\hat
A\otimes_A\cO(X)$ is a finite normal extension. It follows that $\hat
B=\hat B_1\times\ldots\times\hat B_s$ where each $\hat B_i\cong
E_i\[[t^{1/n_i}]]$ with $n_i\in\ZZ_{>0}$ and $[E_i:E]<\infty$.  In
that case, we have that $\hat\cA=\hat A\otimes_A\cA=\cA\otimes_A\hat
A$ is a graded cofinite subalgebra of $\cD_c(\hat B)=\cD_c(\hat
B_1)\times\ldots\times\cD_c(\hat B_s)$ with base $\hat A$. Finally, we
may choose $\tilde Y\pfeil Y$ in such a way that $E_i=E$ for all $i$:
let $D$ be the image of $\tilde D$ in $Y$. Assume the preimage of $D$
in $X$ has irreducible components $D_1,\ldots,D_r$. Then it suffices
to require that $E=\CC(\tilde D)$ is a splitting field for all the
finite extensions $\CC(D_j)|\CC(D)$.

\beginsection relative. Certain rings of differential operators

In this section, we are going to construct a certain class of graded
cofinite subalgebra (see \cite{gradedcofinite}). Later we show that,
under mild conditions, all examples are basically of this kind
(\cite{main}).

For a dominant morphism $\phi:X\pfeil Y$ we have $\cO(Y)\into\cO(X)$
and we can define the subalgebra
$$
\cD(X,Y):=\{D\in\cD(X)\mid D(\cO(Y))\subseteq\cO(Y)\}.
$$
Its associated graded algebra is denoted by $\cDq(X,Y)$.

Now assume that the field extension $\CC(X)|\CC(Y)$ is finite. This
means that there is a non\_empty open subset $X_0\subseteq X$ such
that $\phi:X_0\pfeil Y$ is \'etale. Then every differential operator
$D$ on $Y$ can be uniquely lifted to differential operator $D_0$ on
$X_0$ (see \cite{E7}). Thus, $\cD(X,Y)$ can be also interpreted as the set of
$D\in\cD(Y)$ such that $D_0$ extends to a (regular) differential
operator on $X$. In other words, the diagram
$$
\matrix{\cD(X,Y)&\into&\cD(X)\cr
\uinto&&\uinto\cr
\cD(Y)&\into&\cD(X_0)\cr}
$$
is cartesian. Note that the filtrations of $\cD(X,Y)$ induced by
those on $\cD(X)$, $\cD(Y)$, and $\cD(X_0)$ are the same. Thus we get
an analogous diagram of inclusions for the associated graded rings
$$
\matrix{\cDq(X,Y)&\into&\cDq(X)\cr
\uinto&&\uinto\cr
\cDq(Y)&\into&\cDq(X_0)\cr}
$$
which may not be cartesian, however.

First we show that this class of algebras includes rings of invariant
differential operators:

\Theorem finite. Let $W$ be a finite group acting on $X$. Then
$$
\cD(X,X/W)=\cD(X)^W\hbox{ and }\cDq(X,X/W)=\cDq(X)^W.
$$

\Proof: Clearly $\cD(X,X/W)\supseteq\cD(X)^W$. Conversely, for
$D\in\cD(X,X/W)$ put $D'={1\over|W|}\sum_{w\in W}wD$. If
$f\in\cO(X)^W$ then $(wD)(f)=w(D(w^{-1}f))=f$. This implies that
$D-D'$ is a differential operator which is zero on $\cO(X)^W$,
hence on all of $\cO(X)$. Thus, $D=D'\in\cD(X)^W$.  The second
equality follows from the fact that forming the associated graded
algebra commutes with taking $W$\_invariants.\qed

\Example: Let $X=\A^1$ be the affine line with coordinate ring
$\cO(X)=\CC[x]$ and $W=\mu_n\cong\ZZ/n\ZZ$ acting by
multiplication. Define $Y\cong\A^1$ by $\cO(Y)=\CC[t]$ where
$t=x^n$. The chain rule yields
$\partial_x=nx^{n-1}\partial_t=nt^{1-{1\over n}}\partial_t$. Let
$\xi$ and $\tau$ be the symbols of $\partial_x$ and $\partial_t$,
respectively. Then $\zeta\in W$ acts on $(x,\xi)$ by
$(\zeta^{-1}x,\zeta\xi)$. Moreover, $\xi=nx^{n-1}\tau$. Thus, we have
$$
\vbox{\halign{&\hfill$#$\ \hfill\cr
\cD(X,Y)&=&\CC\<t,t\partial_t,(t^{1-{1\over n}}\partial_t)^n\>&=&\CC\<x^n,x\partial_x,\partial_x^n\>&\into&\CC\<x,\partial_x\>&=&\cD(X)\cr
\uinto&&\uinto&&\uinto&&\uinto&&\uinto\cr
\cD(Y)&=&\CC\<t,\partial_t\>&=&\CC\<x^n,x^{1-n}\partial_x\>&\into&\CC\<x,x^{-1},\partial_x\>&=&\cD(X_0)\cr}}
$$
while the corresponding diagram for the associated graded rings is
$$6
\vbox{\halign{&\hfill$#$\ \hfill\cr
\cDq(X,Y)&=&\CC[t,t\tau,t^{n-1}\tau^n]&=&\CC[x^n,x\xi,\xi^n]&\into&\CC[x,\xi]&=&\cDq(X)\cr
\uinto&&\uinto&&\uinto&&\uinto&&\uinto\cr
\cDq(Y)&=&\CC[t,\tau]&=&\CC[x^n,x^{1-n}\xi]&\into&\CC[x,x^{-1},\xi]&=&\cDq(X_0)\cr}}
$$

\medskip In general, not all subalgebras of the form
$\cD(X,Y)$ are graded cofinite. To formulate a criterion we introduce
the following notions.

\Definition: Let $X$ and $Y$ be normal varieties and $\phi:X\pfeil Y$
a finite surjective morphism. Let $D\subseteq Y$ be a prime divisor
and consider the divisor $\phi^{-1}(D)=r_1E_1+\ldots+r_sE_s$ where the
$E_i$ are pairwise distinct prime divisors and $r_i>0$. We say that
$\phi$ is {\it uniformly ramified over $D$} if $r_1=\ldots=r_s$. Moreover,
$\phi$ is uniformly ramified if it is uniformly ramified over every
$D\subseteq Y$. If all the ramification numbers $r_i$ are $1$ for all
$D$ then we call $\phi$ {\it unramified in codimension one}. Equivalently,
there is an open subset $X_0\subseteq X$ with
$\|codim|_X(X\setminus X_0)\ge2$ on which $\phi$ is \'etale.

\Proposition unramified. Let $X\pfeil Y$ be a finite dominant morphism
between normal affine varieties which is unramified in codimension one. Then
$\cD(X,Y)=\cD(Y)$.

\Proof: Let $D\in\cD(Y)$. Then $D$ can be uniquely lifted to a
differential operator $D_0$ on the set $X_0\subseteq X$ on which
$\phi$ is \'etale. Since $\|codim|_X(X\setminus X_0)\ge2$ and since
$X$ is normal we have $\cO(X_0)=\cO(X)$. Hence one can extend $D_0$
uniquely to all of $X$ which proves $\cD(Y)\subseteq\cD(X,Y)$.\qed

\noindent
Now we show that uniformly ramified morphisms are just a slight
generalization of quotients by finite groups. For this we introduce
the following notation: let $W$ be a finite group acting on a normal
variety $X$. For a prime divisor $Z\subseteq X$ let $W_Z\subseteq W$
be the pointwise stabilizer of $Z$ in $W$ (the inertia group). This
group is always a cyclic group and its order is the ramification number
of $X\pfeil X/W$.

Now assume that $\phi:X\pfeil Y$ is a finite surjective morphism between
normal varieties. Then the field extension $\CC(X)|\CC(Y)$ is finite,
hence has a Galois cover $L$ with Galois group $W$. Let $H\subseteq W$
be the Galois group of $L|\CC(X)$ and let $\tilde X$ be
the normal affine variety such that $\cO(\tilde X)$ is the integral
closure of $\cO(Y)$ in $L$. Then $\tilde X$ carries a $W$\_action with
$\tilde X/W=Y$ and $\tilde X/H=X$ and we have the diagram
$$
\matrix{\tilde X&\pfeil&\tilde X/H&=&X\cr
&\searrow&\downarrow&&\Links\phi\downarrow\cr
&&\tilde X/W&=&Y\cr}
$$

\Proposition reduction. Using the notation above, the
following statements are equivalent:
\Item{8}The morphism $\phi:X\pfeil Y$ is uniformly ramified.
\Item{9}The morphism $\tilde X\pfeil\tilde X/H$ is unramified in codimension
one.
\Item{10} For all prime divisors $Z$ of $\tilde X$ the condition
$W_Z\cap H=1$ holds.\Par
\noindent Moreover, under these conditions holds $\cD(X)=\cD(\tilde X)^H$ and
$\cD(X,Y)=\cD(\tilde X)^W$.

\Proof: The inertia group of $\tilde X\pfeil\tilde X/H$ at $Z$ is
$H_Z=W_Z\cap H$ which shows the equivalence
\cite{I9}$\Leftrightarrow$\cite{I10}.

Let $D$ be the image of $Z$ in $Y$. Then the divisors of $\tilde X$
lying over $D$ are precisely the translates $wZ$, $w\in W$. For fixed
$w\in W$ let $E$ be the image of $wZ$ in $\tilde X/H$. Then $E$ is a
prime divisor of $X$ lying over $D$ and every such divisor is of this
kind.

The inertia group of $\tilde X\pfeil\tilde X/H$ and $\tilde X\pfeil
X/W$at $wZ$ is $H\cap W_{wZ}=H\cap wW_Zw^{-1}$ and $W_{wZ}=wW_Zw^{-1}$
respectively. Therefore, the ramification number of $X\pfeil Y$ at $E$
is $[wW_Zw^{-1}:H\cap wW_Zw^{-1}]$. Thus, condition \cite{I8} means that the
order of $w^{-1}Hw\cap W_Z\cong H\cap wW_Zw^{-1}$ is independent of
$w\in W$. This means in turn that all isotropy groups of $W_Z$ acting
on $W/H$ have the same order. Now $W_Z$, being cyclic, has at most one
subgroup of any given order. Therefore, \cite{I8} means that all isotropy
groups in $W_Z$ on $W/H$ are the same.  We assumed that $L|\CC(Y)$ is
the Galois cover of $\CC(X)|\CC(Y)$, i.e., the smallest Galois
extension of $\CC(Y)$ containing $\CC(X)$. This means that $H$ does
not contain any non\_trivial normal subgroup of $W$, i.e., that the
action of $W$ on $W/H$ is effective. We conclude that \cite{I8} is equivalent
to the statement that for all $Z$ the isotropy groups of $W_Z$ on
$W/H$ are trivial. This is precisely the content of \cite{I9}.

Finally, $\cD(X)=\cD(\tilde X,X)=\cD(\tilde X)^H$ follows from
\cite{unramified} and \cite{finite}. Moreover,
$D\in\cD(X,Y)$ implies $D\in\cD(X)=\cD(\tilde
X)^H\subseteq\cD(\tilde X)$. Hence $D\in\cD(\tilde X,Y)=\cD(\tilde
X)^W$. This shows $\cD(X,Y)\subseteq\cD(\tilde X)^W$. The opposite
inclusion is obvious.\qed

\noindent This result makes it easy to construct uniformly ramified
morphisms. Take, for example, $W=S^n$ with its standard action on the
affine space $\A^n$. Let $H\subseteq S_n$ be a subgroup of odd order
(or any other subgroup not containing a transposition). Then the
morphism $\A^n/H\pfeil\A^n/S_n$ is uniformly ramified. For any {\it
fixed} $X$ it appears to be quite difficult to construct uniformly
unramified morphisms. For $X=\A^1$, the affine line, see \cite{Aeins}
and its proof.

The following technical consequence will be crucial in the
proof of \cite{simple}.

\Corollary retraction. For any uniformly ramified morphism
$\phi:X\pfeil Y$ the inclusion $\cD(X,Y)\into\cD(X)$ has a left inverse
$\rho:\cD(X)\auf\cD(X,Y)$ which is a homomorphism of $\cD(X,Y)$\_bimodules.

\Proof: The map $\rho$ is just the averaging operator
${1\over|W|}\sum_ww$ restricted
to $H$\_invariants $\cD(\tilde X)^H\pfeil\cD(\tilde X)^W$.
\qed

Now we show that uniform ramification is also necessary for
$\cD(X,Y)$ to be graded cofinite. In fact, we prove something
stronger:

\Theorem uniframif. Let $\phi:X\pfeil Y$ be a dominant morphism
between normal affine varieties. Assume $\cD(X)$ contains a graded cofinite
subalgebra $\cA$ with base $\cO(Y)$. Then $\phi$ is uniformly
ramified.

\Proof: That $\phi$ is finite follows from \cite{finitebase}. Suppose
that $\phi$ is not uniformly ramified and let $D\subseteq Y$ be a
prime divisor over which $\phi$ has non\_uniform ramification. Choose
$\tilde Y\pfeil Y$ \'etale and a prime divisor $\tilde D\subset\tilde
Y$ which maps to $D$ as in the example at the end of
section~\cite{basechange}. Then $\hat\cA=\hat A\otimes_A\cA$ is graded
cofinite in $\cD_c(\hat B)=\cD_c(\hat B_1)\times\ldots\times\cD_c(\hat
B_s)$. This means that $\hat\cA$ is actually graded cofinite in each
of the algebras $\cD_c(\hat B_i)$.

We have $\hat A=E\[[t]]$ and $\hat B_i=E\[[t^{1/n_i}]]$. Let
$u_2,\ldots,u_m$ be a transcendence basis of $E$ and put $u_1:=t$. Let
$\partial_i$ be the associated partial derivatives of $\hat A$. Their
symbols are denoted by $\eta_i$ with the special notation
$\tau:=\eta_1$. Then $\cDq_c(\hat A)=E\[[t]][\tau,S]$ with
$S=\{\eta_2,\ldots,\eta_m\}$. Similarly to the example after
\cite{finite} we have $\cDq_c(\hat B_i,\hat
A)=E\[[t]][t\tau,t^{n_i-1}\tau^{n_i},S]$ (see diagram~\cite{E6}).

Since $\phi$ is non\_uniformly ramified over $D$ the $n_i$ are not all
equal. Hence, after relabeling we may assume $n_1<n_2$. Then we see
that $\cDq_c(\hat B_2,\hat A)\subseteq \cDq_c(\hat B_1,\hat A)$
(considered as subrings of $\cDq_c(\hat A)$). Since
$\overline{\hat\cA}\subseteq
\cDq_c(\hat B_1,\hat A)\cap\cDq_c(\hat B_2,\hat A)$ is cofinite in
$\cDq_c(\hat B)$ we conclude that $\cDq_c(\hat B_2,\hat A)$ is
cofinite in $\cDq_c(\hat B_1)$. Now put $x:=t^{1/n_1}$ and let $\xi$
be the symbol of $\partial_x$. Then $\xi=n_1x^{n_1-1}\tau$ implies
$$
n_1^{n_2}t^{n_2-1}\tau^{n_2}=x^{(n_2-1)n_1+n_2(1-n_1)}\xi^{n_2}=x^{n_2-n_1}\xi^{n_2}
$$
and therefore
$$
\cDq_c(\hat B_2,\hat A)=E\[[t]][t\tau,t^{n_2-1}\tau^{n_2},S]=
E\[[x^{n_2}]][x\xi,x^{n_2-n_1}\xi^{n_2},S]\subseteq
E\[[x]][\xi,S]=\cDq_c(\hat B_1).
$$
Now put $x=0$. Then $\cDq_c(\hat B_2,\hat A)$ becomes $E[S]$ which is
clearly not cofinite in $E[\xi,S]$.\qed

\Definition: The affine variety $X$ is called {\it $\cD$\_finite\/} if
$\cDq(X)$ is a finitely generated $\CC$\_algebra.

\medskip \noindent All smooth varieties are $\cD$\_finite. The cubic
$x^3+y^3+z^3=0$ is the standard example of a variety which is not
$\cD$\_finite (see \cite{BGG}).

\Corollary gradedcofinite. Let $\phi:X\pfeil Y$ be a dominant morphism between
normal affine varieties and assume $X$ to be $\cD$\_finite. Then the
following are equivalent:
\Item{} $\phi$ is uniformly ramified.
\Item{} $\cD(X,Y)$ is graded cofinite in $\cD(X)$.\Par

\Proof: If $\phi$ is uniformly ramified then
$\cDq(X)\subseteq\cDq(\tilde X)$ is integral over $\cDq(X,Y)=\cDq(\tilde
X)^W$ (notation as in \cite{reduction}). Since $\cDq(X)$ is finitely
generated it is even finite over $\cDq(X,Y)$. The converse follows
from \cite{uniframif}.\qed

\beginsection simplicity. Simplicity

In this section we derive a simplicity criterion for the ring
$\cD(X,Y)$.

\Definition: A $\cD$\_finite affine variety $X$ is called {\it
$\cD$\_simple} if $\cD(X)$ is a simple ring.

\medskip\noindent It is well known that smooth varieties are
$\cD$\_simple. A curve is $\cD$\_simple if and only if its
normalization map is bijective (see \cite{SmSt}). Further examples
include quotients $X/W$ of smooth varieties by finite groups. More
generally, Schwarz conjectured, \cite{Sch}, that any categorical
quotient $X\mod G:=\|Spec|\cO(X)^G$ is $\cD$\_finite where $G$ is a
reductive group and $X$ is a smooth $G$\_variety. This has been
confirmed in many cases (\cite{Sch1}, \cite{Sch2}, \cite{VdB2}). It
should be added that $\cD$\_simple varieties are automatically
Cohen\_Macaulay (Van den Bergh \cite{VdB1}). In particular, for a
$\cD$\_simple variety normality is equivalent to smoothness in
codimension one.

\Lemma inter. For an affine variety $Y$ let $I\subseteq\cD(Y)$ be a
non\_zero subspace with $[\cO(Y),I]\subseteq I$. Then
$I\cap\cO(Y)\ne0$.

\Proof: Let $0\ne D\in I$ be of minimal order. From minimality and
$[\cO(Y),D]\in I$ we get $[\cO(Y),D]=0$, i.e.,
$D\in\|End|_{\cO(Y)}\cO(Y)=\cO(Y)$.\qed

\Theorem simple. Let $X\pfeil Y$ be a uniformly ramified morphism
between normal affine varieties.
\Item{5} $X$ is $\cD$\_finite if and only if $\cDq(X,Y)$ is finitely
generated.
\Item{6} If $X$ is $\cD$\_simple then $\cD(X,Y)$ is simple.
\Item{7} If $\cD(X,Y)$ is simple then $\cD(X)$ and $\cD(Y)$ are simple.

\Proof: \cite{I5} The algebra $S:=\cDq(X)$ is integral over
$R:=\cDq(X,Y)$. Thus, if $S$ is finitely generated then $R$
is so as well by \cite{finitecomm}. Let $L$ be the field of
fractions of $S$. Then $L|\CC$ is a finitely generated field
extension. Thus, if $R$ is finitely generated then its integral
closure in $L$ is a finite $R$\_module. This implies that $S$ is
finitely generated.

\cite{I6} Let $I\subseteq\cD(X,Y)$ be a non\_zero two\_sided ideal. By
\cite{inter} we may choose a non\_zero function $f\in I\cap\cO(Y)$.
\cite{gradedcofinite}, states that $\cDq(X)$ is a finitely generated
$\cDq(X,Y)$\_module. This implies that $\cD(X)$ is a finitely
generated right $\cD(X,Y)$\_module. In other words, there are
operators $D_1,\ldots,D_s\in\cD(X)$ such that
$\cD(X)=\sum_iD_i\cD(X,Y)$. Let $k$ be an integer which is strictly
larger than the order of every $D_i$. Then
$(\|ad|f)^k(D_i)=0$. On the other side we have
$(\|ad|f)^k(D_i)=\sum_{\nu=0}^k(-1)^\nu{k\choose\nu}f^{k-\nu}D_if^\nu$,
hence $f^kD_i\in\cD(X)f\subseteq\cD(X)I$. This means that $f^k$
annihilates the $\cD(X)$\_module $\cD(X)/\cD(X)I$. The annihilator is
a two\_sided ideal and $\cD(X)$ is a simple ring, thus
$\cD(X)=\cD(X)I$. Applying the retraction $\cD(X)\auf\cD(X,Y)$ from
\cite{retraction} shows $\cD(X,Y)=\cD(X,Y)I=I$.

\cite{I7} Let $Z$ denote either $X$ or $Y$ and assume that $I\ne0$ is
a two\_sided ideal of $\cD(Z)$. \cite{inter} implies that there is
$0\ne f\in I\cap\cO(Z)$. Since $Z\pfeil Y$ is finite we have
$L:=\cO(Z)f\cap\cO(Y)\ne0$. Since $L\subseteq I\cap\cD(X,Y)$ this
shows $I\cap\cD(X,Y)\ne0$. Hence $1\in I\cap\cD(X,Y)\subseteq I$.\qed

\Corollary. Let $X\pfeil Y$ be a finite morphism between normal affine
varieties which is unramified in codimension one. Then
$X$ is $\cD$\_simple if and only if $Y$ is.

\Proof: In this case is $\cD(X,Y)=\cD(Y)$.\qed

\noindent Thus, if one wishes then one may assume in the following
that we are always in the situation that $X\pfeil Y$ is a quotient by
a finite group.

\def\oO{{\overline\Omega}}

\beginsection Poisson. The associated graded algebra

The next result is the beginning of the classification of all graded
cofinite subalgebras.

\Lemma gradedfinite. For a field extension $L|\CC$ let
$\cAq\subseteq L[\xi_1,\ldots,\xi_n]$ be a cofinite homogeneous
subalgebra. Then its base $K=L\cap\cAq$ is a field. Moreover:
\Item3 Assume ${\partial\over\partial\xi_i}\cAq\subseteq\cAq$ for all
$i=1,\ldots,n$. Then $\cAq=K[\xi_1,\ldots,\xi_n]$.
\Item4 Assume ${\partial\over\partial\xi_i}\cAq\subseteq\cAq$ just
for $i=1,\ldots,n-1$. Then there is a positive integer $k$ and
$a_1,\ldots,a_n\in L$ such that
$K[\xi_1+a_1\xi_n,\ldots,\xi_{n-1}+a_{n-1}\xi_n,a_n\xi_n^k]
\subseteq\cAq$.\Par

\Proof: Clearly, $K$ is cofinite in $L$. This implies that
$IL\subseteq L$ is a proper ideal whenever $I\subset K$ is a proper
ideal. This forces $I=0$ and implies that $K$ is a field.

\cite{I3} For a multiindex $\alpha\in\NN^n$ define
$\xi^\alpha=\xi_1^{\alpha_1}\ldots\xi_n^{\alpha_n}$ and analogously
$\partial^\alpha$. Let $f=\sum_\alpha c_\alpha\xi^\alpha\in\cAq$ be
homogeneous. Then $\partial^\alpha(f)=\alpha!c_\alpha\in L\cap\cAq=K$ which
implies $\cAq\subseteq K[\xi,\ldots,\xi_n]$.

For the reverse inclusion it suffices to show $\xi_i\in\cAq$ for all
$i$. Let $S$ be the intersection of $\cAq$ with
$\<\xi_1,\ldots,\xi_n\>_K$. Then, after a linear change of coordinates,
we may assume $S=\<\xi_1,\ldots,\xi_m\>_K$ for some $m\le n$. Since
$\cAq$ is cofinite there is a homogeneous $f\in\cAq$ such that the
variable $\xi_n$ occurs in $f$. Assume the monomial $\xi^\alpha$
occurs in $f$ with $\alpha_n>0$. Put $\beta=\alpha-e_n$ where
$e_n=(0,\ldots,0,1)$. Then $\partial^\beta(f)$ an element of $S$ which
contains $\xi_n$. This implies $m\ge n$ and we are done.

\cite{I4} Let $\pi:L[\xi_1,\ldots,\xi_n]\pfeil
L[\xi_1,\ldots,\xi_{n-1}]$ be the projection obtained by setting
$\xi_n=0$. Then part \cite{I3} implies that $\pi(\cAq)$ contains
$\xi_1,\ldots,\xi_{n-1}$. Thus $\cAq$ contains elements of the form
$\xi_i+a_i\xi_n$ for $i=1,\ldots,n-1$.

Now perform the coordinate change $\xi_i\mapsto\xi_i-a_i\xi_n$ for
$i=1,\ldots,n-1$ and $\xi_n\mapsto\xi_n$. This is allowed since the
partial derivatives $\partial/\partial\xi_i$, $i=1,\ldots,n-1$ stay
unchanged. So we may assume $a_1=\ldots=a_{n-1}=0$.

Since $\cAq$ is cofinite there is an element $f=\sum_\alpha
c_\alpha\xi^\alpha$ which contains the variable $\xi_n$, i.e.,
$c_\alpha\ne0$ and $\alpha_n>0$ for some multiindex $\alpha$. Assume
that $k:=\alpha_n>0$ is as small as possible. Put
$\beta:=(\alpha_1,\ldots,\alpha_{n-1},0)$. Then
$g=\partial^\beta(f)\in\cAq$ is of the form
$g=a_n\xi^k+h(\xi_1,\ldots,\xi_{n-1})$. Moreover, each coefficient of
$h$ appears as a derivative $\partial^\gamma(g)$ for a convenient
multiindex $\gamma$ with $\gamma_n=0$. This implies $h\in
K[\xi_1,\ldots,\xi_{n-1}]\subseteq\cAq$ and we are done.\qed

\beginsection cofinite. Automorphisms

For a normal affine variety $X$
let $\Omega(X)$ be the module of K\"ahler differentials. Then
$$
\cT(X)=\|Hom|_{\cO(X)}(\Omega(X),\cO(X))
$$
is the module of vector fields and we have a canonical isomorphism
$$
\cD(X)_{\le1}=\cO(X)\oplus\cT(X).
$$
Let $\oO(X):=\|Hom|_{\cO(X)}(\cT(X),\cO(X))$, the double dual of
$\Omega(X)$. Since $X$ is normal, elements of $\oO(X)$ can be
characterized as those rational $1$\_forms on $X$ which are regular in
codimension one (or, equivalently, on the smooth part $X^s$ of
$X$). Let $Z(X)$ be the set of $\omega\in\oO(X)$ with
$d\omega|_{X^s}=0$. Our interest in $Z(X)$ comes from the following
well\_known

\Lemma. For every $\omega\in Z(X)$ there is a unique automorphism
$\Phi_\omega$ of $\cD(X)$ with $\Phi_\omega(f)=f$ for all $f\in\cO(X)$
and $\Phi_\omega(\xi)=\xi+\omega(\xi)$ for
all $\xi\in\cT(X)$. This automorphism induces the identity on $\cDq(X)$.

\Proof: First assume $X$ to be smooth. Then $\cD(X)$ is generated by
$\cO(X)\cup\cT(X)$ subject to the relations
$$
\xi f-f\xi=\xi(f)\quad\hbox{and}\quad\xi\eta-\eta\xi=[\xi,\eta].
$$
The first relation is clearly satisfied by $\Phi_\omega$. The second
relation is preserved because of Cartan's formula
$$
0=d\omega(\xi,\eta)=\omega([\xi,\eta])-\xi(\omega(\eta))+\eta(\omega(\xi)).
$$
This shows that $\Phi_\omega$ exists. Clearly, it is the identity on
$\cDq(X)$.

In general, we have shown that $\Phi_\omega(D)$ is a differential
operator on the smooth part of $X$. By normality, it is regular on all
of $X$ and still induces the identity on $\cDq(X)$.\qed

\noindent Let $\phi:X\pfeil Y$ be a finite morphism of normal
varieties. Subsequently, we want to study the twists
$\cA=\Phi_\omega\cD(X,Y)$ with $\omega\in Z(X)$. Clearly, $\cA$
doesn't determine $\omega$ since $\Phi_\omega\cD(X,Y)=\cD(X,Y)$ if
$\omega\in Z(Y)$. To pin down a unique $\omega$ we consider the trace
map $\|tr|_{L|K}:L\pfeil K$ where $K$ and $L$ are the function fields
$\CC(Y)$ and $\CC(X)$. This map induces a trace maps
$\Omega(L)\pfeil\Omega(K)$ characterized by the property
$$
\|tr|_{L|K}(f\phi^*\omega)=\|tr|_{L|K}(f)\omega,\quad f\in
L,\omega\in\Omega(K).
$$
It commutes with the derivative $d$ and splits, up to the factor
$[L:K]$, the inclusion $\Omega(K)\into\Omega(L)$. We define
$Z_K(L)$ as the set of $\omega\in\Omega(L)$ with $d\omega=0$ and
$\|tr|_{L|K}\omega=0$.

Recall the following property of the trace: let $\partial_K:K\pfeil K$
be a derivation and $\partial_L:L\pfeil L$ its unique extension to
$L$. Then
$$14
\|tr|_{L|K}(\partial_Lf)=\partial_K\|tr|_{L|K}(f),\quad f\in L.
$$
Indeed, we may assume that $L|K$ is Galois with group $\Gamma$. Then
$\|tr|_{L|K}f=\sum_{\gamma\in\Gamma}\gamma(f)$. Since the extension
$\partial_L$ is unique, it commutes with $\Gamma$ and the claim follows.

All notions have global counterparts: there are induced trace
maps $\cO(X)\pfeil\cO(Y)$ and $\oO(X)\pfeil\oO(Y)$ (see
\cite{Za}). We put $Z_Y(X)=Z_K(L)\cap\oO(X)$.

In the next result, we are classifying graded cofinite subalgebras of
$\cD(X)$ generically:

\Proposition generic. Let $X$ be an affine variety with quotient field
$\CC(X)=L$ and let $\cA\subseteq\cD(L)$ be a graded cofinite algebra
with base $K=\cA\cap L$. Then $K$ is a field with
$[L:K]<\infty$. Furthermore, there is a unique $\omega\in Z_K(L)$
with $\cA=\Phi_\omega\cD(K)$.

\Proof: That $K\subseteq L$ is a cofinite subfield is proved in the
same way as in \cite{gradedfinite}. Let $u_1,\ldots,u_n\in K$ be a
transcendence basis. Then there are unique derivations
$\partial_1,\ldots,\partial_n$ of $L$ (or $K$) with
$\partial_i(u_j)=\delta_{ij}$. Moreover, these derivations together
with $L$ generate the ring $\cD(L)$. Let $\xi_i$ be the symbol of
$\partial_i$. Then we have
$$
K\subseteq \cAq\subseteq\cDq(L)=L[\xi_1,\ldots,\xi_n].
$$
Observe that $\cDq(L)$ is a Poisson algebra and $\cAq$ is a
sub-Poisson algebra. We have $\{f,u_i\}={\partial
f\over\partial\xi_i}$ which means that $\cAq$ is stable under the
operators $\partial/\partial\xi_i$. \cite{gradedfinite} implies
$\cAq=K[\xi_1,\ldots,\xi_n]$. This means in particular that $\cA$
contains elements of the form $\delta_i:=\partial_i+b_i$ with $b_i\in
L$. We may replace $b_i$ by the unique element of $b_i+K$ with trace
zero. If $\omega:=b_1\,du_1+\ldots+b_n\,du_n$ then
$\|tr|_{L|K}\omega=0$.

Observe
$a_{ij}:=[\delta_i,\delta_j]=\partial_i(b_j)-\partial_j(b_i)\in
\cA\cap L=K$. From $\|tr|_{L|K}a_{ij}=0$ (see \cite{E14}) we infer
$a_{ij}=0$. This means $d\omega=0$ and therefore $\omega\in
Z_K(L)$. Since $\delta_i=\Phi_\omega(\partial_i)$ we get
$\Phi_\omega\cD(K)\subseteq\cA$. From $\cAq=\cDq(K)$ we get
$\cA=\Phi_\omega\cD(K)$.\qed

The 1-form $\omega$ from \cite{generic} may have poles. Our goal is to show
that this won't happen if it comes from a graded cofinite subalgebra $\cA$
of $\cD(X)$. First, a very local version of this result:

\Lemma extension. Let $E$ be a finitely generated field extension of
$\CC$ and put $B=E\[[x]]$ and $A=E\[[t]]\subseteq B$ with $t=x^p$ for
some integer $p\ge1$. Let $K=E\((t))$ and $L=E\((x))$ be the fields of
fractions of $A$ and $B$. For $\omega\in Z(L)$ assume that
$\|tr|_{L|K}\omega$ is regular at $t=0$ and that 
$\cA=\cD_c(B)\cap\Phi_\omega\cD_c(K)$ is graded cofinite in
$\cD_c(B)$. Then $\omega$ is regular at $x=0$.

\Proof: If $p=1$ then $\omega=\|tr|_{L|K}\omega$ is regular. Assume
$p\ge2$ from now on. Let $u_1,\ldots,u_{n-1}$ be a transcendence basis
of $E$ and put $u_n=x$. Let $\partial_1,\ldots,\partial_n$ be the
corresponding differentials of $B$. Put $b_i:=\omega(\partial_i)\in
L$. Then we have to show $b_i\in B$ for all $i$. We have
$$
\omega=\sum_{i=1}^{n-1}b_idu_i+b_ndt^{1/p}=
\sum_{i=1}^{n-1}b_idu_i+\textstyle{1\over p}b_nxt^{-1}dt.
$$
Hence the condition that $\|tr|_{L|K}\omega$ is regular means
$$8
\|tr|_{L|K}b_1,\ldots,\|tr|_{L|K}b_{n-1},t^{-1}\|tr|_{L|K}xb_n\in A.
$$

Note also the explicit formula
$$
\|tr|_{L|K}x^d=\cases{px^d=pt^{d/p}&if $p|d$\cr0&otherwise\cr}
$$

Let $\xi_i\in\cDq_c(B)$ be the symbol of $\partial_i$. Then
$\cAq\subseteq\cDq_c(B)=B[\xi_1,\ldots,\xi_n]$ is cofinite. Since
$A\subseteq\cA$ we have $u_1,\ldots,u_{n-1}\in\cAq$. As in the proof
of \cite{generic} this implies that $\cAq$ is stable under partial
differentiation by $\xi_i$, $i=1,\ldots,n-1$. Let $\cAq'$ be the image
of $\cAq$ in $B/xB=E[\xi_1,\ldots,\xi_n]$. Then
\cite{gradedfinite}\cite{I4} applied to $\cAq'$ gives elements
$a_{ij}\in\delta_{ij}+xB$, $i,j=1,\ldots,n-1$ and $c_i\in B$,
$i=1,\ldots,n-1$ such that
$$
\sum_{j=1}^{n-1}a_{ij}\xi_j+c_i\xi_n\in\cAq\quad{\rm for\ }i=1,\ldots,n-1.
$$
From $\cA\subseteq\Phi_\omega\cD_c(K)$ we infer (since
$\partial_n=px^{p-1}\partial_t$)
$$10
\cAq\subseteq\cDq_c(K)=K[\xi_1,\ldots,\xi_{n-1},x^{1-p}\xi_n].
$$
This implies in particular $a_{ij}\in K\cap B=A$. Since
$A\subseteq\cAq$ and since the matrix $(a_{ij})\in M_{n-1}(A)$ is
invertible we may assume $a_{ij}=\delta_{ij}$, i.e.,
$$
\xi_1+c_1\xi_n,\ldots,\xi_{n-1}+c_{n-1}\xi_n\in\cAq\quad
{\rm with}\quad c_1,\ldots,c_{n-1}\in x^{1-p}K\cap B=xA.
$$
In the last equation we used $p\ge2$. Lifting to $\cA$ we get operators
$$12
\delta_i:=\partial_i+c_i\partial_n+d_i\in\cA\quad
{\rm with}\quad c_i\in xA,\ d_i\in B,\ i=1,\ldots,n-1.
$$
Now we use $\Phi_{-\omega}(\cA)\subseteq\cD_c(K)$. More precisely, from
$$
\Phi_{-\omega}(\delta_i)=(\partial_i-b_i)+c_i(\partial_n-b_n)+d_i
$$
we get $d_i-b_i-c_ib_n\in K$. Therefore, $b_i\in K+B+A(xb_n)$. From
\cite{E8} we obtain
$$9
b_i\in B+Axb_n,\quad i=1,\ldots,n-1.
$$

Now we use that \cite{gradedfinite}\cite{I4} gives us also an element
of $\cAq$ of the form $a\xi_n^k\|mod|x$ with $a\in B^\Times$. From
\cite{E10} we obtain
$$
a\in(x^{1-p})^kK\cap B
$$
Since $a$ has a non\_zero constant term this is only possible if
$p$ divides $k$. Then $a\in K\cap B^\Times=A^\Times$, hence we can make
$a=1$. Summarizing, we have found an operator $D$ in $\cA$ of the form
$$13
D=\partial_n^k+(u_1\partial_1+\ldots+
u_{n-1}\partial_{n-1}+u)\partial_n^{k-1}+\ldots\quad
{\rm with}\ u_i\in xB,u\in B\ {\rm and}\ p|k.
$$
As above we want to use that $\Phi_{-\omega}(D)\in\cD_c(K)$. More
precisely we want to look at the coefficient of
$\partial_t^{k-1}$. Write $\partial_n=f\partial_t$ with
$f=px^{p-1}=pt^{1-{1\over p}}$. Using the easily verified formulas
$$
\partial_n^k=(f\partial_t)^k=f^k\partial_t^k+\alpha_{k,p}f^{k-1}x^{-1}\partial_t^{k-1}+\ldots\quad{\rm
with}\ \alpha_{k,p}=(p-1)\textstyle{k\choose2}
$$
$$
\partial_x^{k-1}=f^{k-1}\partial_t^{k-1}+\ldots
$$
$$
(\partial_n-b_n)^k=
f^k\partial_t^k+f^{k-1}(\alpha_{k,p}x^{-1}-kb_n)\partial_t^{k-1}+\ldots
$$
the coefficient of $\partial_t^{k-1}$ in $\Phi_{-\omega}(D)$ can be computed:
$$
f^{k-1}(\alpha_{k,p}x^{-1}-kb_n-u_1b_1-\ldots-u_{n-1}b_{n-1}+u)\in K
$$
From \cite{E9} we get elements $v\in B$ with $w\in B^\Times$ such that
$$
f^{k-1}(\alpha_{k,p}x^{-1}-wb_n+v)\in K
$$
Since $(p-1)(k-1)\equiv1\|mod|p$ we have $x/f^{k-1}\in K$. This implies
$$11
wxb_n\in K+B=K+xB=E\((x^p))+xE\[[x]]\quad.
$$
Let $d\in\ZZ$ be the order of zero of $xb_n$ or, equivalently,
$wxb_n$. If $d\le0$ then \cite{E11} implies $p|d$. On the other hand,
$\|tr|(xb_n)\in tE\[[t]]$ (see \cite{E8}) means that $xb_n$ doesn't
contain any monomials $x^d$ with $p|d$ and $d\le0$. Therefore $d>0$,
i.e., $b_n\in B$. Finally, \cite{E9} implies that the other $b_i$ are
in $B$ and we are done.\qed

\noindent The next statement is similar but much easier to prove:

\Lemma cofinite2. Let $B=E\[[x]]$ with quotient field $L=E\((x))$ and
let $\omega\in Z(L)$. Assume that $\cA=\cD_c(B)\cap\Phi_\omega\cD_c(B)$ is
cofinite in $\cD_c(B)$. Then $\omega\in Z(B)$.

\Proof: The base of $\cA$ is $E\[[x]]$. Thus we have
$u_1,\ldots,u_{n-1},u_n=x\in\cA$ and we can apply right away part
\cite{I3} of \cite{gradedfinite}. Thus we get $\cAq\ \|mod|
x=E[\xi_1,\ldots,\xi_n]$. The Nakayama lemma implies
$\cAq=\cDq_c(E\[[x]])$, hence $\cA=\cD_c(E\[[x]])$. In particular
$\Phi_\omega(\partial_i)=\partial_i+\omega(\partial_i)\in\cD(E\[[x]])$
means that $\omega$ is regular.\qed

\noindent Now we globalize these local computations:

\Theorem regular. Let $\phi:X\pfeil Y$ be a finite dominant morphism
between normal varieties and let $L$, $K$ be the fields of rational
functions of $X$, $Y$ respectively. For $\omega\in Z_K(L)$ assume that
$\cA=\cD(X)\cap\Phi_\omega\cD(K)$ is graded cofinite in $\cD(X)$. Then
$\omega\in Z_Y(X)$, i.e., $\omega$ is regular on all of $X$.

\Proof: Since $X$ is normal it suffices to prove the regularity of
$\omega$ in codimension one. Let $D\subseteq Y$ be a prime divisor and
choose $\tilde D\subseteq\tilde Y\pfeil Y$ as in the final remark of
section~\cite{basechange}. \cite{uniframif} implies that $\phi$ is
uniformly ramified. Therefore, the rings $\hat B_i$ are all the same,
say equal to $E\[[x]]$ with $x^p=t$. The form $\omega$ gives rise to
forms $\omega_i$ over $\hat B_i[x^{-1}]\cong E\((x))$. From
$\cA\subseteq\cD(X)\cap\Phi_\omega\cD(K)$ (with $K=\CC(Y)$) we get
$$
\hat\cA\subseteq\cD_c(E\[[x]])^s
\cap(\Phi_{\omega_1}\times\ldots\times\Phi_{\omega_s})\Delta\cD_c(E\((t)))
$$
where $\Delta$ is the diagonal embedding. Thus $\hat\cA$ is contained in
the set of all $(D_1,\ldots,D_s)\in\cD_c(E\[[x]])^s$ with
$$
\Phi_{-\omega_1}(D_1)=\ldots=\Phi_{-\omega_s}(D_s)\in\cD_c(E\((t))).
$$
Solving for $D_2,\ldots,D_s$ we see that $\hat\cA$ is contained in
$$
\cD_c(E\[[x]])\cap\Phi_{\omega_2-\omega_1}\cD_c(E\[[x]])\cap\ldots\cap
\Phi_{\omega_s-\omega_1}\cD_c(E\[[x]]).
$$
In particular, the latter algebra is graded cofinite in $\cD(X)$ which
implies, by \cite{cofinite2}, that $\delta_i:=\omega_i-\omega_1$ is
regular for all $i$. Then
$$
0=\|tr|_{L|K}\omega=s\|tr|_{E\((x))|E\((t))}\omega_1+
\sum_i\|tr|_{E\((x))|E\((t))}\delta_i.
$$
implies that $\|tr|_{E\((x))|E\((t))}\omega_1$ is regular. Since clearly
$\hat\cA\subseteq\cD_c(E\[[x]])\cap\Phi_{\omega_1}\cD_c(E\((t)))$ we
deduce from \cite{extension} that $\omega_1$ itself is regular.\qed

\beginsection applications. The main theorem and its applications

The main result of this paper is:

\Theorem main. Let $X$ be a normal $\cD$\_simple variety and let
$\cA\subseteq\cD(X)$ be a graded cofinite subalgebra. Then
$\cA=\Phi_\omega\cD(X,Y)$ where $Y=\|Spec|\cA\cap\cO(X)$ and
$\omega\in Z_Y(X)$ unique. The variety $Y$ is normal and the morphism
$X\pfeil Y$ is uniformly ramified.

\Proof: Let $\cA'$ be the normalization of $\cA$ with base $A'$
(\cite{normalization}) . Put $Y':=\|Spec|A'$ and let $L$, $K$ be the
quotient fields of $X$, $Y'$, respectively. \cite{uniframif} implies
that $X\pfeil Y'$ is uniformly ramified. We conclude that
$\cD(X,Y')$ is simple (\cite{simple}).

From \cite{generic} we get
a unique $\omega\in Z_L(K)$ such that
$\cA'_K=\Phi_{\omega}\cD(K)$. By \cite{regular}, this $\omega$ is
regular on all of $X$ and we may replace $\cA$ by
$\Phi_{-\omega}\cA$. Thereby, we get
$$2
\cA\subseteq\cD(X)\cap\cD(K)=\cD(X,Y')\subseteq\cD(X)
$$
We have $K\otimes_A\cA=\cD(K)=K\otimes_A\cD(X,Y')$. Hence, for
every $D\in\cD(X,Y')$ there is $0\ne f\in A$ such that
$fD\in\cA$. Now \cite{E2} implies that $\cD(X,Y')$ is a
finitely generated $\cA$\_module, both left and right. Thus
there is a single $0\ne f\in A$ with
$f\cD(X,Y')\subseteq\cA$. Likewise, there is $0\ne g\in A$ with
$\cD(X,Y')g\subseteq\cA$. This implies that
$\cD(X,Y')gf\cD(X,Y')$ is a non\_zero two\_sided ideal of $\cD(X,Y')$
which is contained in $\cA$. We conclude
$\cA=\cD(X,Y')$. From this we get
$\cO(Y)=\cO(X)\cap\cA=\cO(Y')$, hence $Y=Y'$ is normal.\qed

For the applications we start with a well\_known cofiniteness criterion:

\Lemma cofincrit. Let $R=\oplus_{d=0}^\infty R_d$ be a finitely
generated graded $\CC$\_algebra. Let $F\subseteq R_0$ be a subset such
that $R_0$ is finite over $\CC[F]$. Let $G\subseteq R_{>0}$ be a set
of homogeneous elements which has the same zero\_set in $\|Spec|R$ as
$R_{>0}$. Then the subalgebra generated by $F\cup G$ is cofinite in
$R$.

\Proof: Hilbert's Nullstellensatz implies that there is $N>0$ with
$(R_{>0})^N\subseteq RG$. Since $R$ is finitely generated there is an
$M\ge N$ with $R_{>M}\subseteq(R_{>0})^N$. Put $S:=R_{\le M}$. This is
a finitely generated $R_0$\_module with $R=S+RG$. Thus we have
$$
R=S+RG=S+SG+RG^2=\ldots=S+SG+\ldots+SG^{d-1}+RG^d
$$
for all $d\ge1$. Since the minimal degree of an element of $G^d$ goes
to $\infty$ as $d$ goes to $\infty$ we see that $R=S[G]$, hence is a finitely
generated $R_0[G]$\_module. Thus it is also a finitely generated
$\CC[F\cup G]$\_module.\qed

\Theorem mainappl. Let $W$ be a finite group acting on the normal $\cD$\_simple
affine variety $X$. Let $F\subseteq\cO(X)^W$ and $G\subseteq\cD(X)^W$
with
\Item1 The normalization of $\CC[F]$ is $\cO(X)^W$.
\Item2 The set of symbols $\Gq$ of $G$ vanishes simultaneously only on
the zero section of the cotangent bundle $\|Spec|\cDq(X)$ of
$X$.\Par\noindent
Then $\cD(X)^W$ is, as an algebra, generated by $F\cup G$.

\Proof: Let $\cA\subseteq\cD(X)$ be the subalgebra generated by $F$
and $G$. Then $F$ and $\Gq$ meet the assumptions of \cite{cofincrit}
and we conclude that $\cA$ is graded cofinite in $\cD(X)$.

Let $A$ be the base of $\cA$. By \cite{main} it is integrally
closed. We have $\cA\subseteq\cD(X)^W$ hence $\CC[F]\subseteq
A\subseteq\cO(X)^W$ which implies $A=\cO(X)^W$. Finally, \cite{main}
implies $\cAq=\cDq(X,X/W)=\cDq(X)^W$. From $\cA\subseteq\cD(X)^W$ we
get $\cA=\cD(X)^W$.\qed

As mentioned in the introduction, we obtain the following result of
Levasseur\_Stafford \cite{LS} as an application:

\Corollary. Let $V$ be a finite dimensional representation of
$W$. Then $\cD(V)^W$ is generated by the invariant polynomials along
with the invariant constant coefficient differential operators.

\noindent Observe that even for vector spaces, \cite{mainappl} is more
general than the Levasseur\_Stafford theorem: it suffices to take
invariant functions which generate the ring of invariants only up to
normalization and invariant constant coefficient operators which
generate all invariant constant coefficient operators up to integral
closure. In practice, this leads to much smaller generating sets. For
example, we get

\Corollary. Let $V$ be an $n$\_dimensional representation of $W$. Then
$\cD(V)^W$ can be generated by $2n+1$ elements.

\Proof: First, choose homogeneous systems of parameters
$f_1,\ldots,f_n$ and $d_1,\ldots,d_n$ of $\cO(V)^W$ and $\cO(V^*)^W$,
respectively. Then choose a generator $f_0\in\cO(V)^W$ of the finite field
extension $\CC(V)^W/\CC(f_1,\ldots,f_n)$. Then $F=\{f_0,\ldots,f_n\}$
and $G=\{d_1,\ldots,d_n\}$ satisfy the assumptions of
\cite{mainappl}.\qed

We need the following

\Definition: Let $\cA\subseteq\cD(X)$ be a graded cofinite
subalgebra with base $A$ and $Y=\|Spec|A$. Then $\cA$ is called {\it
untwisted} if $\cA=\cD(X,Y)$.

\Proposition untwisted. Let $X$ be a normal $\cD$\_simple affine variety, and
$\cA\subseteq\cA'\subseteq\cD(X)$ graded cofinite subalgebras. If
$\cA$ is untwisted then so is $\cA'$.

\Proof: Let $K$, $K'$, and $L$ be the field of fractions of
$\cA\cap\cO(X)$, $\cA'\cap\cO(X)$, and $\cO(X)$,
respectively. Moreover, let $\cA_K$ (resp. $\cA'_{K'}$) be the algebra
generated by $\cA$ and $K$ (resp. $\cA'$ and $K'$). Choose a
transcendence basis $u_1,\ldots,u_n\in K$ and let
$\partial_1,\ldots,\partial_n$ be the derivations of $K$, $K'$, and $L$ with
$\partial_i(u_j)=\delta_{ij}$. If $\cA$ is untwisted then
$\partial_i\in\cA_K$. Thus $\partial_i\in\cA'_{K'}$ which means that
$\cA'$ is untwisted, as well.\qed

Now we derive a Galois correspondence for graded cofinite
subalgebras:

\Theorem Galois. Let $X$ be a normal $\cD$\_simple affine variety and $W$ a
finite group acting on $X$. Then the map $H\mapsto\cD(X)^H$
establishes a bijective correspondence between subgroups of $W$ and
subalgebras of $\cD(X)$ containing $\cD(X)^W$.

\Proof: The only non\_trivial thing to show is that every subalgebra
$\cA$ containing $\cD(X)^W$ is of the form $\cD(X)^H$. Let $A$ be the
base of $\cA$ and $Y=\|Spec|A$. By \cite{main}
and \cite{untwisted} we have $\cA=\cD(X,Y)$. Since $\cO(X)^W\subseteq
A\subseteq\cO(X)$ and since $A$ is integrally closed there is
$H\subseteq W$ with $A=\cO(X)^H$. Thus $\cA=\cD(X,X/H)=\cD(X)^H$.\qed

\Remark: The preceding result could have been as well derived from a
noncommutative version of Galois theory due to Kharchenko. Recall that
a subalgebra $\cA$ of $\cD(X)$ is called an {\it anti-ideal} if for
any $a\in\cD(X)$, $b,c\in\cA\setminus\{0\}$, $ab,ca\in\cA$ implies
$a\in\cA$ (see e.g. \cite{Cohn} \S6.6, p.334). This is a
non\_commutative version of integral closedness. Now \cite{Galois}
follows from Kharchenko's Galois correspondence (\cite{Cohn}
Thm.~11.7) using the following

\Proposition. Let $X$ be a normal $\cD$\_simple affine variety. Then
every graded cofinite subalgebra of $\cD(X)$ is an anti-ideal.

\Proof: Let $\cA\subseteq\cD(X)$ be a graded cofinite subalgebra. By
\cite{main} we may assume that $\cA=\cD(X,Y)$ for some uniformly
ramified morphism $X\pfeil Y$. Now using \cite{reduction} we have
$$
\cA=\cD(\tilde X)^W\subseteq\cD(X)\subseteq\cD(\tilde X).
$$
Is is easy to see (see the proof of \cite{Cohn}
Thm.~11.7) that $\cD(\tilde X)^W$ is an anti-ideal of $\cD(\tilde
X)$. A fortiori, it is an anti-ideal of $\cD(X)$.\qed

Here is another example of how one can play with \cite{main}:

\Theorem. Let $V$ be a finite dimensional representation of
$W$. Then the ring $\cD(V\oplus V^*)^W$ is generated by
$$
\cD(V\oplus0)^W\cup\cD(0\oplus V^*)^W\cup\{\omega\}.
$$
where $\omega:V\times V^*\pfeil\CC$ is the evaluation map.

\Proof: Let $\cA$ be the subalgebra generated by this set. The first
two pieces generate the subalgebra $\cD(V\oplus V^*)^{W\times
W}$. \cite{Galois} implies that $\cA=\cD(V\oplus V^*)^H$ for some subgroup
$H$ of $W\times W$. But the isotropy group of $\omega$ inside $W\times
W$ is just $W$ embedded diagonally which implies $H=W$.\qed

One remarkable feature of subalgebras of non\_commutative
rings is that they are much scarcer. An argument similar to
\cite{Galois} shows

\Corollary. Let $X$ be normal and $\cD$\_simple and $\cA\subseteq\cD(X)$ graded
cofinite. Then there are only finitely many intermediate subalgebras.

\Proof: By applying an automorphism to $\cD(X)$ we may assume $\cA$ to
be untwisted: $\cA=\cD(X,Y)$. Then every intermediate algebra is
untwisted as well hence of the form $\cD(X,Y')$ with
$\cO(Y)\subseteq\cO(Y')\subseteq\cO(X)$ and $\cO(Y')$ integrally
closed. Galois theory tells us that there are only finitely
many of those.\qed

For $X=\A^1$ one can make things very explicit:

\Theorem Aeins. Let $\cA\subseteq\cD(\A^1)=\CC\<x,\partial_x\>$ be
graded cofinite. Then there is $a\in\CC$, $p\in\CC[x]$, and
$m\in\ZZ_{>0}$ such that $\cA=\CC\<u^m,\eta^m\>$ where $u=x-a$
and $\eta=\partial_x+p(x)$. Moreover, $p$ may be chosen in such a way
that $xp$ does not contain monomials whose exponent is divisible by
$m$. In that case, the triple $(a,m,p)$ is uniquely determined by $\cA$.

\Proof: Clearly we may assume $\cA$ to be untwisted. Then we have to
determine all uniformly ramified morphism $\phi:\A^1\pfeil Y$.

First, $Y$ is a smooth rational curve with $\cO(Y)^\Times=\CC^*$ which
implies $Y\cong\A^1$. Thus, $\phi$ extends to a morphism
$\phi:\P^1\pfeil\P^1$ such that $\phi^{-1}(\infty)=\infty$. Let $d$
be the degree of $\phi$. Assume $\phi$ is ramified over the points
$y_1,\ldots,y_s\in Y$ with ramification numbers $r_1,\ldots,r_s\ge2$. Then
$\phi^{-1}(y_i)$ will consist of $d/r_i$ points. The ramification
number at $\infty$ is $d$. Thus Hurwitz' formula implies
$$
-2=-2d+\sum_{i=1}^s{d\over\ r_i}(r_i-1)+(d-1)=(s-1)d-1-\sum_{i=1}^s{d\over r_i}
$$
From $r_i\ge 2$ we get
$$
-1=(s-1)d-\sum_{i=1}^s{d\over r_i}\ge(s-1)d-s{d\over2}=({s\over2}-1)d.
$$
This implies $s=0$ and $d=1$, i.e., $\phi$ is an isomorphism, or $s=1$
and $r_1=d$. In the latter case, $\phi$ is, up to a translation, just
the quotient $\A^1\pfeil\A^1/\mu_d$.\qed

\noindent {\bf Final remark:} The stipulation that our subalgebras are
{\it graded\/} cofinite in $\cD(X)$ is essential. It would be
interesting to classify all subalgebras $\cA$ for which $\cD(X)$
itself is a finitely generated left and right $\cA$\_module. Take, for
example, the affine space $X=\A^n$. Then $\cD(X)$ is the Weyl algebra
on which the symplectic group $Sp_{2n}(\CC)$ acts by
automorphisms. Now take any irreducible $2n$\_dimensional
representation of a finite group $W$ which preserves a symplectic
form. Then $\cA=\cD(\A^n)^W$ will have the required property even
though it is not graded cofinite. The point is, of course, that the
$W$\_action does not preserve the standard filtration. Nevertheless,
it preserves the so\_called Bernstein filtration for which linear
functions have degree one. Therefore, one might want to start with the
problem: what are the subalgebras of a Weyl algebra which are graded
cofinite with respect to the Bernstein filtration?

\beginrefs

\L|Abk:BGG|Sig:BGG|Au:Bernstein, J.; Gelfand, I.; Gelfand,
S.|Tit:Differential operators on a cubic cone|Zs:Uspehi Mat.
Nauk|Bd:27|S:185--190|J:1972|xxx:-||

\B|Abk:Cohn|Sig:Co|Au:Cohn, P.|Tit:Free rings and their
relations. Second edition.|Reihe:London Mathematical Society
Monographs {\bf 19}|Verlag:Academic Press, Inc. [Harcourt Brace
Jovanovich, Publishers]|Ort:London|J:1985|xxx:-||

\B|Abk:Eisen|Sig:Ei|Au:Eisenbud, David|Tit:Commutative
algebra. With a view toward algebraic geometry|Reihe:Graduate Texts in
Mathematics {\bf 150}|Verlag:Springer-Verlag|Ort:New York|J:1995|xxx:-||

\L|Abk:LS|Sig:LS|Au:Levasseur, T.; Stafford, J.|Tit:Invariant
differential operators and an homomorphism of Harish-Chandra|Zs:J.
Amer. Math. Soc.|Bd:8|S:365--372|J:1995|xxx:-||

\L|Abk:Masson|Sig:Mas|Au:M\'asson, G.|Tit:Rings of differential
operators and \'etale homomorphisms|Zs:MIT Thesis|Bd:-|S:-|J:1991|xxx:homepage.mac.com/gisli.masson/thesis/||

\L|Abk:Sch1|Sig:\\Sch|Au:Schwarz, G.|Tit:Differential operators on
quotients of simple groups|Zs:J. Algebra|Bd:169|S:248--273|J:1994|xxx:-||

\L|Abk:Sch2|Sig:\\Sch|Au:Schwarz, G.|Tit:Lifting differential operators
from orbit spaces|Zs:Ann. Sci. \'Ecole Norm. Sup.
(4)|Bd:28|S:253--305|J:1995|xxx:-||

\Pr|Abk:Sch|Sig:\\Sch|Au:Schwarz, G.|Artikel:Invariant differential
operators|Titel:Proceedings of the International Congress of
Mathematicians (Z\"urich, 1994)|Hgr:-|Reihe:-|Bd:-%
|Verlag:Birkh\"auser|Ort:Basel|S:Vol. {\bf1}, 333--341|J:1995|xxx:-||

\L|Abk:SmSt|Sig:SS|Au:Smith, S.; Stafford, J.|Tit:Differential operators on an affine curve|Zs:Proc. London Math. Soc.|Bd:56|S:229--259|J:1988|xxx:-||

\Pr|Abk:VdB1|Sig:\\VdB|Au:Van den Bergh, M.|Artikel:Differential operators
on semi-invariants for tori and weighted projective spaces|Titel:Topics in
Invariant Theory|Hgr:M.-P. Malliavin ed.|Reihe:Lecture Notes in
Math.|Bd:1478|Verlag:Springer|Ort:Berlin|S:255--272|J:1991|xxx:-||

\L|Abk:VdB2|Sig:\\VdB|Au:Van den Bergh, M.|Tit:Some rings of differential
operators for ${\rm Sl}_2$-invariants are simple. {\rm Contact Franco-Belge
en Alg\`ebre (Diepenbeek, 1993)}|Zs:J. Pure Appl.
Algebra|Bd:107|S:309--335|J:1996|xxx:-||

\L|Abk:Wa|Sig:Wa|Au:Wallach, N.|Tit:Invariant differential operators
on a reductive Lie algebra and Weyl group representations|Zs:J. Amer.
Math. Soc.|Bd:6|S:779--816|J:1993|xxx:-||

\L|Abk:Za|Sig:Za|Au:Zannier, U.|Tit:A note on traces of differential
forms|Zs:J. Pure Appl. Algebra|Bd:142|S:91--97|J:1999|xxx:-||

\endrefs

\bye